\documentclass[10pt,draft]{book}

\usepackage{amsmath, amssymb}
\usepackage{amsfonts, mathrsfs}
\usepackage[cp1251]{inputenc}
\usepackage[active]{srcltx}
\usepackage[ukrainian, russian, english]{babel}

\begin{document}

\selectlanguage{ukrainian} \setcounter{page}{214}

\makeatletter
\renewcommand{\@evenhead}{\footnotesize{\tt
Збірник\,праць\,Ін-ту\,математики\,НАН\,України \hfill
2014,\,Т.11,\,№ 3,\,\thepage--\pageref{end}}}
\renewcommand{\@oddhead}{\footnotesize{\tt
Збірник\,праць\,Ін-ту\,математики\,НАН\,України \hfill
2014,\,Т.11,\,№ 3,\,\thepage--\pageref{end}}} \makeatother

\noindent {\small УДК 517.5} \vskip 1.5mm

\noindent \textbf{В.\,В.~Савчук} {\small (Ін--т математики НАН України, Київ)} \vskip 1.5mm

\noindent \textbf{{МНОГОЧЛЕНИ ФАБЕРА ЗІ СПІЛЬНИМ КОРЕНЕМ}}

\vskip 2mm

{\small \it \noindent
We describe the two sets of meromorphic univalent functions in the class $\Sigma$, for which the sequence of  Faber polynomials $\{F_j\}_{j=1}^\infty $ have the roots with following properties respectively: $\sum_{j=1}^{n}|F_j(z_0)|=0<|F_{n+1}(z_0)|, $ $n\in\mathbb N, $ and $|F_1(z_0)|> 0=\sum_{j=2}^{\infty}|F_j(z_0)|$. We found an explicit form of Faber polynomials for such functions.
\hfill}

\vskip 1.5mm
{\small \it \noindent Описано дві множини мероморфних однолистих функцій класу $\Sigma$, для яких послідовності многочленів Фабера $\{F_j\}_{j=1}^\infty$ мають корені відповідно з такими властивостями : $\sum_{j=1}^{n}|F_j(z_0)|=0<|F_{n+1}(z_0)|,$ $n\in\mathbb N,$ і $|F_1(z_0)|>0=\sum_{j=2}^{\infty}|F_j(z_0)|$. Знайдено явний вигляд многочленів Фабера для таких функцій.
\hfill}

 \normalsize \vskip 3mm

{\bf Вступ.}
Нехай $\Sigma$ --- клас функцій
$$
\Psi(w)=w+\sum_{j=0}^\infty\alpha_jw^{-j},\eqno(1)
$$
які є мероморфними і однолистими в області $\mathbb D^-:=\{w\in\widehat{\mathbb C} : |w|>1\}$.

Многочленами Фабера функції $\Psi\in\Sigma$ називається послідовність алгебраїчних многочленів $\{F_j\}_{j=1}^\infty$, які визначаються як коефіцієнти розкладу
$$
\ln\frac{\Psi(w)-z}{w}=-\sum_{j=1}^\infty \frac{F_j(z)}{j}w^{-j}, \quad z\in\mathbb C,\eqno(2)
$$
в степеневий ряд відносно $w$  в околі нескінченно віддаленої точки (див., наприклад, [1, p. 57]).

Многочлени Фабера можна означити і для довільної мероморфної функції $\Psi$, не обов'язково однолистої. Таке означення дається за допомогою рекурентної формули (див., наприклад, [2, с. 60]): {\it системою многочленів Фабера мероморфної в $\mathbb D^-$ функції $\Psi$, яка має розклад (1), називається система $\mathcal F(\Psi):=\{F_j\}_{j=0}^\infty$ алгебраїчних мно\-гочленів $F_j$ степеня $j$ $(F_0(z)=1, F_1(z)=z-\alpha_0)$ таких, що для
}
\vskip 3mm \hfill \textbf{\copyright\ В.\,В.~Савчук, 2014}

\noindent
{\it будь-якого $z\in\mathbb C$ і кожного натурального $j$ справджуються
рівності
$$
F_{j+1}(z)+(\alpha_0-z)F_{j}(z)+\sum_{k=1}^{j}\alpha_k F_{j-k}(z)+j\alpha_j=0.\eqno(3)
$$
}

Зауважимо, що з рекурентних рівностей (3) для даної системи многочленів $\{F_j\}_{j=0}^\infty$ і послідовності комплексних чисел $\{\alpha_j\}_{j=0}^\infty$, взагалі кажучи, не випливає той факт, що послідовність $\{\alpha_j\}_{j=0}^\infty$  є послідовністю коефіцієнтів Лорана--Тейлора мероморфної однолистої функції $\Psi$. Якщо ж розглядати означення многочленів Фабера, дані на основі співвідношень (2) і (3) на класі $\Sigma$, то вони рівносильні. Зазначимо також, якщо існує система алгебраїчних многочленів така, що для заданої послідовності комплексних чисел  $\{\alpha_j\}_{j=0}^\infty$  виконується (3), то така система єдина.

В численних застосуваннях многочленів Фабера, зокрема, в теорії наближення аналітичних функцій комплексної змінної, важливою є інформація про корені цих многочленів. Дослідженням питань про розміщення, асимптотичний розподіл коренів многочленів Фабера тощо, присвячено багато робіт. Бібліографію з цього кола задач можна знайти в [3].

Мета даної роботи --- описання множини функцій $\Psi\in\Sigma$, для яких тільки перші $n$  многочленів Фабера мають один спільний корінь а також описання множини функцій, для яких усі многочлени Фабера, починаючи з деякого, і тільки вони мають один спільний корінь.

Робота написана за такою схемою.

В п. 1 дано опис множини функцій $\Psi\in\Sigma$, для яких $\sum_{j=1}^{n}|F_j(z_0)|=0<|F_{n+1}(z_0)|$. Також описано функції $\Psi$, для яких відрізок послідовності многочленів Фабера $\{F_j\}_{j=0}^m$ збігається з відрізком послідовності многочленів Тейлора $\{(z-\alpha_0)^k\}_{k=0}^{m}$, а відрізок $\{F_j\}_{j=m+1}^n,$ $m<n-1$, з відрізком послідовності алгебраїчних многочленів, які задовольняють тричленне рекурентне співвідношення зі сталим коефіцієнтом.

В п. 2 показано, що єдиною функцією $\Psi\in\Sigma$, для якої послідовність $\mathcal F(\Psi)$ має спільний корінь $z_0$ починаючи з многочлена $F_2$  є функція $\Psi(w)=z_0+w\exp((\alpha_0-z_0)/w)$.

В п. 3 знайдено явний вигляд многочленів Фабера функції $\Psi(w)=z_0+w\exp((\alpha_0-z_0)/w)$, а також встановлено деякі їх властивості.

\bigskip

{\bf 1.} Добре відомо, що для функції $\Psi(w)=w+\alpha_0$ многочлени Фабера мають вигляд $F_j(z)=(z-\alpha_0)^j,$ $j=0,1,2,\ldots$. Тому точка $z_0=\alpha_0$ є спільним коренем для всіх $F_j$ , $j=1,2,\ldots$. В наступному твердженні розвинуто це спостереження.

{\bf Теорема 1.} {\it Нехай $\Psi\in\Sigma$, $\mathcal F(\Psi)=\{F_j\}_{j=0}^\infty$ --- система многочленів Фабера функції $\Psi$, $n\in\mathbb N\cup\{\infty\}$ і $z_0\in\mathbb C$. Тоді наступні твердження рівносильні:

1) $\displaystyle\sum_{j=1}^{n}|F_j(z_0)|=0<|F_{n+1}(z_0)|;$

2) $\Psi(w)=w+z_0+\displaystyle\sum_{j=n}^\infty\alpha_jw^{-j}\quad\forall~w\in\mathbb D^{-}, |\alpha_n|>0~\bigg(\sum_{j=\infty}^\infty=0\bigg);$

3) $F_j(z)=(z-z_0)^j,$ $j=\overline{0, n}$ і $F_{n+1}(z)\not\equiv(z-z_0)^{n+1}~\forall$ $z\in\mathbb C.$
}

{\bf\textsl {Зауваження} 1.} Коефіцієнти Лорана--Тейлора функції $\Psi\in\Sigma$, про яку йдеться в твердженні 2) теореми 1 задовольняють таку нерівність [4, p. 139]:
\[
|\alpha_j|\le\frac{2}{j+1},\quad j=\overline{n, 2n}.
\]
Цікаво зауважити також (це випливає з (3)), що для такої функції~$\Psi$
\[
\alpha_j=-\frac{F_{j+1}(z_0)}{j+1},\quad j=\overline{n, 2n}.
\]

\textbf{\emph{Доведення}.} Твердження теореми 1 є тривіальним при $n=1$, тому далі вважаємо, що $n\ge 2$.

{\it Доведення $"1)\Rightarrow 2)".$} Згідно з (3)
\[
-(j+1)\alpha_j=F_{j+1}(z_0)+(\alpha_0-z_0)F_j(z_0)+\sum_{k=1}^{j-1}\alpha_kF_{j-k}(z_0),\quad j\in\mathbb N,
\]
де (і скрізь далі) суми вигляду $\sum_{j=n}^m$ при $m<n$ покладаються рівними нулю.

\newpage

Права частина останньої рівності дорівнює нулю для всіх натуральних $j\le n-1$, до того ж $F_1(z_0)=z_0-\alpha_0=0$, тому $\alpha_0=z_0$ і $\alpha_j=0,$ $j=\overline{1, n-1}$.

Оскільки $|F_{n+1}(z_0)|>0$, то $|\alpha_n|=(n+1)^{-1}|F_{n+1}(z_0)|>0.$

{\it Доведення $"2)\Rightarrow 3)".$} Згідно з (3)
\[
F_{j}(z_0)=(z-z_0)F_{j-1}(z),\quad j=\overline{1, n}
\]
і
\[
F_{n+1}(z)=(z-z_0)F_n(z)-(n+1)\alpha_n.
\]

Отже,
\[
F_j(z)=(z-z_0)^{j},\quad j=\overline{0, n}
\]
і
\[
F_{n+1}(z)=(z-z_0)^{n+1}-(n+1)\alpha_n\not\equiv (z-z_0)^{n+1}.
\]

{\it Імплікація $"3)\Rightarrow 1)"$} є очевидною.

\bigskip

Повертаючись до теореми 1, зауважимо, що  для функцій $\Psi(w)=w+\alpha_0+\sum_{k=n}^\infty\alpha_kw^{-k}$ і тільки для них перші $n$ многочленів Фабера збігаються з першими $n$ многочленами Тейлора $(z-\alpha_0)^k$, які відповідають функції $w\mapsto w+\alpha_0.$ У зв'язку з цим природно виникає питання про описання функцій $\Psi\in\Sigma$, для яких перші $n$ многочленів Фабера збігаються з іншими добре відомим частинними випадками многочленів Фабера, наприклад, з многочленами Чебишева.

Наступне твердження описує деякі з таких випадків.

{\bf Теорема 2.} {\it Нехай $\Psi\in\Sigma$, $\mathcal F(\Psi)=\{F_j\}_{j=0}^\infty$ --- система многочленів Фабера функції $\Psi$, $z_0\in\mathbb C,$ $m, n\in\mathbb N\cup\{\infty\},$ $m<n-1$. Тоді наступні твердження рівносильні:

1) $\displaystyle \sum_{j=1,\atop
j\not=m+1}^n|F_j(z_0)|=0<|F_{m+1}(z_0)F_{n+1}(z_0)|;$

2) $\Psi(w)=w+z_0+\alpha_mw^{-m}+\displaystyle\sum_{j=n}^\infty\alpha_jw^{-j}$ $\quad\forall~w\in\mathbb D^{-},\quad|\alpha_m\alpha_n|>0;$

3)
\[
F_{j+1}(z)=\left\{
\begin{matrix}
(z-z_0)^{j+1},\hfill j=\overline{0,m-1},\cr\cr
(z-z_0)^{m+1}-(m+1)\alpha_m,\hfill j=m,\cr\cr
(z-z_0)F_{j}(z)-\alpha_mF_{j-m}(z),\hfill j=\overline{m+1, n-1},\cr\cr
(z\!-\!z_0)F_{j}(z)\!-\!\alpha_mF_{j-m}(z)\!-\!\displaystyle\sum_{k=n}^j\alpha_kF_{j-k}(z)\!-\!j\alpha_j,\hfill j\ge n.
\end{matrix}
\right.
\]

}

\textbf{\emph{Доведення}} є цілком аналогічним до доведення теореми 1. Тому окреслимо лише його ключові моменти.

Якщо справджується твердження 1), то за теоремою  1 $\alpha_j=0,$ $j=\overline{1, m-1}$ і $F_j(z)=(z-\alpha_0)^j$ $j=\overline{1,m}$. Далі, згідно з (3) $F_{m+1}(z_0)=-(m+1)\alpha_m\not=0$ і $0=F_{j+1}(z_0)=-\sum_{k=m}^{j-1}\alpha_kF_{j-k}(z_0)-ja_j,$ $j=\overline{m, n-1}$, тобто $\alpha_j=0,$ $j=\overline{m+1,n-1}$.

На підставі цих фактів усі рекурентні формули в твердженні 3) випливають із формули (3) і навпаки.

\bigskip
{\bf\textsl {Приклад} 1.}
Нехай в умовах теореми 2 $\alpha_m=1/m$, тобто
$$
\Psi(w)=w+\frac{1}{mw^m}+\sum_{j=n}^\infty\alpha_jw^{-j},\quad m\in\mathbb N.\eqno(4)
$$
Відомо, що функція $\widetilde\Psi(w):=w+1/(mw^m)$ здійснює конформне відображення області $\mathbb D^{-}$ на зовнішність $m$--гіпоциклоїди з $m+1$ точкою звороту (див., наприклад [5, с. 295]).

Для такої функції $\widetilde\Psi$ в [6] одержано явний вираз для системи многочленів Фабера $\mathcal F(\widetilde\Psi)$:
\[
\widetilde F_j(z)=j\sum_{k=0}^{[j/(m+1)]}\frac{(-1)^k\Gamma(j-mk)}{\Gamma(j-(m+1)k+1)m^kk!}z^{j-(m+1)k},~j\in\mathbb N,
\]
де
\[
\left[\frac{j}{m+1}\right]:=\left\{
\begin{matrix}
\displaystyle\frac{j}{m+1},\hfill &&j=0~\mathop{\rm mod}(m+1),\cr\cr
\displaystyle \frac{j-l}{m+1},\hfill &&j=l~\mathop{\rm mod}(m+1),~l=\overline{1, m},
\end{matrix}
\right.
\]
і $\Gamma(\cdot)$ --- гамма--функція Ейлера, а в роботах [7, 8] показано, що всі нулі многочленів $\widetilde F_j$ лежать на променях, що з'єднують точку $0$ з вершинами $m$--гіпоциклоїди.

Зокрема, при $m=1$ областю значень функції $\Psi$ є вся комплексна площина з розрізом вздовж відрізка $[-2,2]$, тобто $\widetilde\Psi(\mathbb D^{-})=\widehat{\mathbb C}\setminus[-2,2]$. Отже,
\[
\mathcal F(\widetilde\Psi)=\{T_0\}\cup\left\{2T_j\left(\frac{z}{2}\right)\right\}_{j=1}^\infty,
\]
де $T_j$ --- многочлени Чебишева (див., наприклад, [9, p. 56, 57]).

Таким чином згідно з теоремою 2 для функції $\Psi$, визначеної рівністю (4),
$F_j(z)=\widetilde F_j(z),$ $j=\overline{0, n}.$

Многочлени Фабера $F_j$, розглянуті в прикладі 1, при $j=\overline{m+1,n}$ утворюють підклас алгебраїчних многочленів, що породжуються, так званими, тричленними рекурентними співвідношеннями. Такі многочлени відіграють важливу роль в багатьох розділах сучасного аналізу, зокрема,  в асимптотичній теорії ортогональних многочленів і в теорії апроксимацій Ерміта--Паде. Цим та іншим аспектам теорії присвячена робота [10].

\bigskip
{\bf 2.}
Розглянемо тепер задачу про {\it описання множини функцій $\Psi\in\Sigma$, для яких існує принаймні одна точка $z_0\in\mathbb C$ така, що для деякого натурального $n$
$$
\sum_{j=1}^n|F_j(z_0)|>0=\sum_{j=n+1}^\infty|F_j(z_0)|.\eqno(5)
$$
}

Наступне твердження дає розв'язок поставленої задачі у випадку $n=1$.

{\bf Теорема 3.} {\it Нехай $\Psi$ --- мероморфна в $\mathbb D^-$ функція вигляду (1), $\mathcal F(\Psi)=\{F_j\}_{j=0}^\infty$ --- система многочленів Фабера функції $\Psi$ і $z_0\in\mathbb C$. Тоді:

1) для того щоб співвідношення (5) виконувалося при $n=1$ необхідно і достатньо, щоб $z_0\not=\alpha_0$ і
$$
\Psi(w)=z_0+w\exp\left(\displaystyle\frac{\alpha_0-z_0}{w}\right)=w+\alpha_0+\displaystyle\sum_{j=1}^\infty\frac{(\alpha_0-z_0)^{j+1}}{(j+1)!}w^{-j};\eqno(6)
$$

2) функція $\Psi$, задана формулою (6), є однолистою тоді і тільки тоді, коли $|\alpha_0-z_0|\le 1.$
}

{\bf Наслідок 1.} {\it Нехай $\Psi$ --- функція мероморфна в $\mathbb D^-$  і $\mathcal F(\Psi)$ --- система многочленів Фабера функції $\Psi$. Якщо $z_0$ --- точка, для якої виконується співвідношення (5) при $n=1$, то вона єдина.
}

Справді, рівняння $F_2(z)=(z-\alpha_0)^2-(z_0-\alpha_0)^2=0$ має два корені: $z_1=z_0,$ $z_2=2\alpha_0-z_0$. Тому іншою точкою, для якої може виконуватися співвідношення (5) є тільки точка $z_0^*:=2\alpha_0-z_0$. Але в такому випадку для всіх $w\in\mathbb D^-$ мала б виконуватися рівність $\exp((\alpha_0-z_0)/w)=\exp(-(\alpha_0-z_0)/w)$, що неможливо, оскільки $z_0\not=~\alpha_0$.

{\bf\textsl {Зауваження} 2.} Як буде видно з доведення теореми 3, умова $|\alpha_0-z_0|\le 1$ є критерієм того, що функція $\Psi$, задана формулою (6), є зірковою відносно точки $z_0$, тобто, що область $\widehat{\mathbb C}\setminus\overline{\Psi(\mathbb D^-)-z_0}$ є зірковою відносно початку координат.

\textbf{\emph{Доведення}.} {\it 1).}
Нехай $F_2(z_0)=F_3(z_0)=\ldots=0$ і $z_0\not=\alpha_0$. Тоді застосовуючи формулу (3) послідовно для $j=1,2,\ldots$ при $z=z_0$, отримаємо рівності
\[
0+(\alpha_0-z_0)(z_0-\alpha_0)+2\alpha_1=0,
\]
$$
0+(\alpha_0-z_0)\cdot\!0+\alpha_1(z_0-\alpha_0)+3\alpha_2=0,\eqno(7)
$$
\[
0+(\alpha_0-z_0)\cdot\!0+\alpha_1\cdot\!0+\alpha_2(z_0-\alpha_0)\cdot\!0+4\alpha_3=0
\]
і т. д.

Звідси випливає, що
\[
\alpha_1=\frac{(z_0-\alpha_0)^2}{2}, \alpha_2=\frac{(z_0-\alpha_0)^3}{2\!\cdot\!3},\alpha_3=\frac{(z_0-\alpha_0)^4}{2\!\cdot3\!\cdot\!4},\ldots.
\]

\newpage
Отже,
\[
\Psi(w)=w+\alpha_0+\sum_{j=1}^\infty\frac{(z_0-\alpha_0)^{j+1}}{(j+1)!}w^{-j}.
\]

Навпаки, якщо функція $\Psi$ має вигляд (6) і $z_0\not=\alpha_0$, то будуть виконуватися рівності (7), які згідно з формулами (3) й доводять, що $F_2(z_0)=F_3(z_0)=\ldots=0$.

Зауважимо, що якщо ж виконується (5) при $n=1$ і функція $\Psi\in\Sigma$, то висновок про її вигляд можна зробити безпосередньо з рівності (2). Справді, в такому випадку розклад (2) в околі нескінченно віддаленої точки набуває вигляду
$$
\ln\frac{\Psi(w)-z_0}{w}=-F_1(z_0)w^{-1}=(\alpha_0-z_0)w^{-1}.\eqno(8)
$$
Припустимо, що $z_0=\Psi(w_0)$ для деякого $w_0\in\mathbb D^-$. Тоді ліву частину рівності (8) можна переписати у вигляді
\[
\ln\frac{\Psi(w)-z_0}{w}=\ln\frac{\Psi(w)-\Psi(w_0)}{w-w_0}+\ln\left(1-\frac{w_0}{w}\right),
\]
звідки випливає, що при значеннях $w$ з околу $w_0$ абсолютна величина лівої частини рівності (8) за рахунок доданка $\ln(1-w_0/w)$ може бути як завгодно великою, в той час як права частина (8) є обмеженою.

Отримана суперечність доводить, що рівність (8) справджується для всіх $w\in\mathbb D^-$  і не існує такого $w_0\in\mathbb D^-$, щоб $z_0=\Psi(w_0)$.

Отже, неодмінно $z_0\in\widehat{\mathbb C}\setminus\overline{\Psi(\mathbb D^-)}$ і
$$
\Psi(w)=z_0+w\exp\left(\frac{\alpha_0-z_0}{w}\right),\quad w\in\mathbb D^-.
$$

{\it 2).} Нехай функція $\Psi$, задана формулою (6), є однолистою.

Тоді згідно з відомою ознакою однолистості (див. [11]) справджується співвідношення
\[
A(\Psi):=\sup_{w\in\mathbb D^-}\left\{(|w|^2-1)\left|w\frac{\Psi''(w)}{\Psi'(w)}\right|\right\}\le 6.
\]

Оскільки
\[
\Psi'(w)=\frac{w-\lambda}{w}\exp\left(\frac{\lambda}{w}\right),~
\Psi''(w)=\frac{\lambda^2}{w^3}\exp\left(\frac{\lambda}{w}\right),\quad w\in\mathbb D^{-},
\]
де $\lambda:=\alpha_0-z_0$, то
\[
6\ge A(\Psi)\ge(R^2-1)\left|\frac{\lambda^2}{Re^{i\arg\lambda}\left(Re^{i\arg\lambda}-\lambda\right)}\right|=
\]
\[
=(R^2-1)\frac{|\lambda|^2}{R\big|R-|\lambda|\big|}\quad\forall~R> 1,
\]
звідки й випливає, що $|\lambda|\le 1$.

Доведемо тепер достатність умови $|\lambda|\le 1$ для однолистості функції $\Psi(w)=z_0+w\exp(\lambda w^{-1})$.

Зрозуміло, що однолистість функції $\Psi$ є рівносильною однолистості функції $w\mapsto(\Psi(w)-z_0)$. Тому, не втрачаючи загальності, покладаємо $z_0=0$.

Насправді, в нашому випадку про функцію $\Psi$ можна ствердити більше. А саме, функція $\Psi$ є зірковою в $\mathbb D^-$, а отже і однолистою. У цьому легко переконатися за допомогою такого критерія (див., наприклад [1, p. 42]): голоморфна функція $\Psi$ є зірковою в $\mathbb D^-$ тоді і тільки тоді, коли
\[
B(\Psi):=\inf_{w\in\mathbb D^-}\mathop{\rm Re}\left(w\frac{\Psi'(w)}{\Psi(w)}\right)\ge 0.
\]

Маємо
\[
B(\Psi)=\inf_{w\in\mathbb D^-}\mathop{\rm Re}\frac{w-\lambda}{w}=\inf_{w\in\mathbb D^-}\left(1-\left|\frac{\lambda}{w}\right|\right)=1-|\lambda|\ge 0,
\]
що й доводить зірковість і однолистість функції $\Psi$.

\bigskip
{\bf 3.} Наведемо тепер деякі властивості многочленів Фабера функції $\Psi$, заданої формулою (6).

Розпочнемо з твердження, в якому дано явний вигляд многочленів Фабера такої функції.

{\bf Теорема 4.} {\it Нехай $\Psi(w)=\eta+w\exp\left(\frac{\lambda}{w}\right)$, $\eta, \lambda\in\mathbb C$, $|\lambda|\le 1$ і $\mathcal F(\Psi)$ --- система многочленів Фабера функції $\Psi$. Тоді
\[
F_1(z)=z-\eta-\lambda,\quad F_j(z)=j\sum_{k=1}^j(-\lambda)^{j-k}\frac{k^{j-k-1}}{(j-k)!}(z-\eta)^k,~j\ge 2,0^0=1.
\]
}

\textbf{\emph{Доведення}.} Теорема є очевидною у випадку, коли $\lambda=0$, тому далі в доведенні цей випадок виключаємо.

Нехай $\Phi:=\Psi^{-1}$ --- функція, обернена до $\Psi$. За теоремою 2 функція $\Psi$ є однолистою в $\mathbb D^-$, тому $\Phi$ є визначеною і однолистою в області $\Psi(\mathbb D^-)$.

Добре відомо (див., наприклад, [2], с. 53), що многочлени Фабера $F_j$ функції $\Psi$ збігаються з правильною частиною (Тейлорівською частиною) функції $\Phi^j$ в розкладі в ряд Лорана в околі нескінченно віддаленої точки, тобто
$$
(\Phi(z))^j=Q_j(z)+F_j(z),\quad j=1,2,\ldots,\eqno(9)
$$
де $Q_j(z)=O(1/z)$, $|z|\to\infty$.

Знайдемо явний вигляд функції $\Phi$. Для цього при довільному фіксованому $z\in\Psi(\mathbb D^-)$ розглянемо рівняння відносно $w$:
\[
z=\eta+w\exp\left(\frac{\lambda}{w}\right).
\]

Оскільки $\eta\not\in\Psi(\mathbb D^-)$, то останнє рівняння рівносильне такому
\[
-\frac{\lambda}{w}\exp\left(-\frac{\lambda}{w}\right)=-\frac{\lambda}{z-\eta}.
\]
Розв'язки цього рівняння мають вигляд
$$
w=-\frac{\lambda}{W\left(-\displaystyle\frac{\lambda}{z-\eta}\right)},\eqno(10)
$$
де $W$ --- функція Ламберта, яка означається як обернена функція до функції $w\mapsto we^w$, тобто $W(we^w)=w$ (див., наприклад, [12]).

Функція $W$ є багатозначною з розгалуженням в точці $-e^{-1}$ і межею однозначності гілок --- променем $(-\infty, -e^{-1}]$. Але оскільки $-\lambda/(z-\eta)\not\in(-\infty,-e^{-1}]$,  то в правій частині (10) беремо головну гілку $W_0$ функції Ламберта.

Таким чином,
$$
\Phi(z)=-\frac{\lambda}{W_0\left(-\displaystyle\frac{\lambda}{z-\eta}\right)}.\eqno(11)
$$

Оскільки, як показано в [12], для будь-якого цілого $j$
\[
\left(W_0(t)\right)^j=-j\sum_{k=j}^{+\infty}\frac{(-k)^{k-j-1}}{(k-j)!}t^k,\quad |t|<e^{-1},~(0/0=1),
\]
то згідно з (9) і (11)
\[
F_j(z)=(-\lambda)^{j}\left(W_0\left(-\frac{\lambda}{z-\eta}\right)\right)^{-j}-Q_j(z)=
\]
\[
=(-\lambda)^{j}j\sum_{k=-j}^0\frac{(-k)^{k+j-1}}{(k+j)!}\left(-\frac{\lambda}{z-\eta}\right)^k,\quad\left(0^0=1\right),
\]
що й потрібно було довести.

\bigskip

Розглянемо задачу з [13, p. 57] про явний вигляд алгебраїчних многочленів $P_j$, які породжуються твірною функцією
\[
K(z,t):=\frac{A(t)}{1-zg(t)},
\]
де $A(t)=\sum_{k=0}^\infty a_kt^k,$ $|a_0|>0$  і $g(t)=\sum_{k=1}^\infty g_kt^k,$ $|g_1|>0$ --- функції, голоморфні в крузі $\mathbb  D:=\{t\in\mathbb C : |t|<1\}$.

Йдеться про послідовність алгебраїчних многочленів $\{P_j\}_{j=0}^\infty$, для яких справджується рівність
$$
K(z,t)=\sum_{j=0}^\infty P_j(z)t^j\eqno(12)
$$
при деякому фіксованому $z$ для всіх $t$ в околі $0$, а також про визначення області $\Omega\times\mathbb D$ всіх таких допустимих значень $(z,t)$.

{\bf Теорема 5.} {\it Нехай $A(t)=1,$ $g(t)=t\exp(-\lambda t),$ $\lambda\in\mathbb C,$ $|\lambda|\le1.$ Тоді рівність (12) справджується для всіх $(z,t)\in\Omega\times\mathbb D$, де $\Omega$ --- зіркова відносно точки $0$ область з межею $\Gamma:=\{\zeta=e^{i\theta}\exp(\lambda e^{-i\theta}): \theta\in[0, 2\pi]\}$,
\[
P_j=\sum_{k=0}^j\lambda^{j-k}F_k,\quad j=0,1,2,\ldots,
\]
і $F_k$ --- многочлени Фабера функції $\Psi(w)=w\exp\left(\frac{\lambda}{w}\right)$.
}

\textbf{\emph{Доведення}.} За теоремою 3 і зауваженням 2 функція $\Psi$ є зірковою, тобто область $\Omega:=\widehat{\mathbb C}\setminus\overline{\Psi(\mathbb D^-)}$ є зірковою відносно точки $0$. Зрозуміло також, що крива $\Gamma$ є межею області $\Omega$.

Нехай $z\in\Omega$. Тоді продиференціювавши відносно $w$ рівність (2), отримаємо добре відоме співвідношення
$$
\frac{\Psi'(w)w}{\Psi(w)-z}=\sum_{j=0}^\infty F_j(z)w^{-j},\eqno(13)
$$
в якому ряд відносно $w$ збігається рівномірно і абсолютно в області $\mathbb D^-$.

З другого боку для всіх $z\in\Omega$ і $w\in\mathbb D^-$
\[
\frac{\Psi'(w)w}{\Psi(w)-z}=\frac{w-\lambda}{w}K\left(z,\frac{1}{w}\right).
\]

Отже, за правилом множення степеневих рядів для будь-якого $z\in\Omega$ і $t\in\mathbb D\setminus\{0\}$ маємо рівності
\[
K(z,t)=\frac{1}{1-\lambda{t}}\cdot\frac{\Psi'\left(\displaystyle\frac{1}{t}\right)
\displaystyle\frac{1}{t}}{\Psi\left(\displaystyle\frac{1}{t}\right)-z}=
\]
\[
=\left(\sum_{k=0}^\infty\lambda^k t^k\right)\left(\sum_{k=0}^\infty F_k(z)t^k\right)=
\sum_{j=0}^\infty\left(\sum_{k=0}^j\lambda^{j-k}F_k(z)\right)t^j,
\]
що й потрібно було довести.

Область $\Omega$ не можна розширити не змінивши область $\mathbb D$. У цьому легко переконатися зауваживши таке.

Припустимо $z\in\Gamma$, тобто $z=\Psi(w_0)$ для деякого $w_0, |w_0|=1$. Тоді точка $w_0$ є особливою точкою для ряду в правій частині (13), що унеможливлює всі висновки, зроблені вище.

Теорему доведено.

{\bf Теорема 6.} {\it Нехай $\Psi(w)=w\exp\left(\frac{\lambda}{w}\right)$, $\lambda\in\mathbb C, |\lambda|\le 1$ і $\mathcal F(\Psi)$ --- система многочленів Фабера функції $\Psi$. Тоді
$$
zF'_j(z)=j\sum_{k=0}^j\lambda^{j-k}F_k(z),\quad j=0,1,2,\ldots.\eqno(14)
$$
}

\textbf{\emph{Доведення}.} За теоремою 5 для будь-якого $z\in\Omega$ і $w\in\mathbb  D^{-}$
\[
\sum_{j=0}^\infty P_j(z)w^{-j}=\frac{1}{1-\displaystyle\frac{z}{w}\exp\left(-\displaystyle\frac{\lambda}{w}\right)}=\frac{w\exp\left(\displaystyle\frac{\lambda}{w}\right)}{w\exp\left(\displaystyle\frac{\lambda}{w}\right)-z}=
\]
$$
=1+\frac{z}{w\exp\left(\displaystyle\frac{\lambda}{w}\right)-z}.\eqno(15)
$$

З другого боку внаслідок диференціювання (2) відносно $z$ (див. також [2, с. 129]) справджується рівність
\[
\frac{1}{\Psi(w)-z}=
\]
$$
=\frac{1}{w\exp\left(\displaystyle\frac{\lambda}{w}\right)-z}=\sum_{j=1}^\infty\frac{F'_{j}(z)}{j}w^{-j}\quad\forall~z\in\Omega, w\in\mathbb  D^{-}.\eqno(16)
$$

Підставивши (16) в (15) і зрівнявши коефіцієнти при $w^{-j}$ в розкладах в обох частинах, за теоремою 5 отримаємо (14).

\bigskip

\footnotesize
\begin{enumerate}
\Rus

\item\label{Pommerenke}
{\it Pommerenke Chr.} Univalent functions. --- G\"{o}ttingen:
Vandenhoeck $\&$ Ruprecht, 1975. --- 376~p.

\item\label{Suetin}
{\it Суетин  П.К.\/} Ряды по многочленам Фабера. ---  М.: Наука, 1984. --- 336~с.

\item\label{Andrievskii-Blatt}
{\it Andrievskii V. V., Blatt H.--P.} On the Distribution of Zeros of Faber Polynomials // Comp. Meth. Funct. Th. --- 2011. --- {\bf 11}, № 1. --- P.~263~---~282.

\item\label{Duren}
{\it Duren P.} Univalent Functions. --- New York: Springer--Verlag, 2001. --- 383~p.

\item\label{Ivanov-Popov}
{\it Иванов В. И., Попов В. Ю.\/} Конформные отображения и их приложения. --- М: Едиториал УРСС, 2002. --- 324 с.

\item\label{He}
{\it  He M. X.\/} Explicit Representations of Faber Polynomials for
m-Cusped Hypocycloids // J. Approx. Th. --- 1996. --- {\bf 87}. --- P.~137~---~147.

\item\label{Eirmann-Varga}
{\it Eiermann M. Varga R. S.} Zeros and local extreme points of Faber polynomials
associated with hypocycloidal domains // Elec. Trans. Numer. Anal. --- 1993. --- {\bf 1}. --- P.~49~---~71.

\item\label{He_Saff}
{\it  He M. X., Saff E. B.\/} The zeros of Faber polynomials for an $m$-cusped hypocycloid // J. Approx. Th. --- 1994. --- {\bf 78}. --- P.~410~---~432.

\item\label{Gil}
{\it Gil A., Segura J., Temme N. M.} Numerical Methods for Special Functions. --- Philadelphia: SIAM, 2007. --- 417 p.

\item
{\it Aptekarev A. I., Kalyagin V. A., Saff  E. B.} Higher--Order Three--Term Recurrences
and Asymptotics of Multiple Orthogonal Polynomials // Constr. Approx. --- 2009. --- {\bf 30}. --- P.~175~---~223.

\item\label{Avhadiev}
{\it Авхадиев Ф. А.} Об условиях однолистности аналитических функций // Изв. вузов. Матем. --- 1970. --- № 11. --- С.~3 --- 13.

\item\label{W function}
{\it Corless R. M., Gonnet G. H., Hare D. E. G., Jeffrey D. J., Knuth D. E.\/} On the Lambert $W$ function // Adv. Comp. Math. --- 1996. --- {\bf 5}. --- P.~329~---~359.

\item\label{Boas-Buck}
{\it Boas R. P. Jr., Buck R. G.} Polynomial expansions of analytic functions . --- Berlin: Springer--Verlag, 1964. --- 77~p.

\end{enumerate}

\label{end}

\end{document}